\documentclass[12pt]{article}
\usepackage{amssymb,amsfonts,amsthm,amsmath,graphicx,fullpage,multirow,longtable,pst-barcode,caption2,lscape}

\newtheorem{thm}{Theorem}

\theoremstyle{remark}

\theoremstyle{definition}
\newtheorem{dfn}{Definition}
\title{Graphs with Non-unique Decomposition and Their Associated Surfaces}
\author{Weiwen Gu}
\date{}
\begin{document}
\maketitle
\begin{abstract}
The ideal (tagged resp.) triangulation of bounded surface with marked points are associated with skew-symmetric (skew-symmetrizable) exchange matrices. An algorithm is established to decompose the graph associated to such matrix. There are finite many graph with non-unique decomposition. We find all such graphs and their decompositions. In addition, we also find the associated ideal (tagged) triangulations to different decompositions.
\end{abstract}
\section{Introduction}
Triangulation is a useful tool to study the topology of surfaces. Ideal triangulation of bordered surfaces with marked points is of particular interests in cluster algebra. For example, in \cite{FST1}, the authors construct cluster algebra associated to an ideal triangulation.
\begin{dfn}
We associate to each ideal triangulation $T$ the (generalized) signed adjacency matrix $B=B(T)$ that reflects the combinatorics of $T$. The rows and columns of $B(T)$ are naturally labeled by the arcs in $T$. For notational convenience, we arbitrarily label these arcs by the numbers $1,\ldots,n$, so that the rows and columns
of $B(T)$ are numbered from $1$ to $n$ as customary, with the understanding that this
numbering of rows and columns is temporary rather than intrinsic. For an arc
(labeled) $i$, let $\pi_T(i)$ denote (the label of) the arc defined as follows: if there is a
self-folded ideal triangle in $T$ folded along $i$, then $\pi_T(i)$ is its remaining
side (the enclosing loop); if there is no such triangle, set $\pi_T(i)=i$. For each ideal triangle $\triangle$ in T which is not self-folded, define the $n\times n$ integer
matrix $B^{\triangle}=(b^{\triangle}_{ij})$ by settings:
\begin{equation*}
b^{\triangle}_{ij}=
\begin{cases}
1 & \text{if $\triangle$ has sides labeled $\pi_T(i)$ and $\pi_T(j)$}\\
 & \text{with $\pi_T(j)$ following $\pi_T(i)$ in the clockwise order;}\\
-1 & \text{if the same holds, with the counter-clockwise order;}\\
0 & \text{otherwise.}
\end{cases}
\end{equation*}
The matrix $B=B(T)=(b_{ij})$ is then defined by $$B=\sum_{\triangle}B^{\triangle}$$
The sum is taken over all ideal triangles $\triangle$ in $T$ which are not self-folded. The $n\times n$ matrix
$B$ is skew-symmetric, and all its entries $b_{ij}$ are equal to $0,1,-1,2$, or $-2$.
\end{dfn}

A quiver is defined as a finite oriented multi-graph without loops and $2$-cycles.

\begin{dfn}
Let $G$ be a quiver, $B(G)=(b_{ij})$ is the skew-symmetric matrix whose rows and columns are labeled by the
vertices of $G$, and whose entry $b_{ij}$ is equal to the number of edges going from $i$ to $j$
minus the number of edges going from $j$ to $i$.
\end{dfn}

\begin{dfn}
Suppose $B$ is a signed adjacency matrix associated to an ideal triangulation of a bordered surface with marked points $(S,M)$, and $G$ is a quiver. If $B(G)=B$, we say $G$ is the \textit{oriented adjacency graph} associated to $(S,M)$.
\end{dfn}

The notion of \textit{Block decomposition} plays an important role in determining the mutation class of a quiver. It is proved in \cite{FST1} that a \textit{quiver} is \textit{block-decomposable} if and only if it is the associated adjacency graph of an ideal triangulations of a bordered surface with marked points. A quiver is a finite oriented multi-graph without loops and $2$-cycles. In \cite{WG}, we provide an algorithm that determines if a given quiver is block decomposable. In addition, we find all connected decomposable graphs with non-unique block-decomposition.

In \cite{FST}, the authors generalize the property to the graph associated to ideal (tagged) triangulation of bordered surfaces with marked points. A new decomposability called \textit{s-decomposable} is studied. It is proved in the same article that there is a one-to-one correspondence between $s$-decomposable skew-symmetrizable graphs with fixed block decomposition and ideal tagged triangulations of marked bordered surfaces with fixed tuple of conjugate pairs of edges. In \cite{WG2}, we provide a generalized algorithm that determines if a given graph is  $s$-decomposable. In addition, we find that only two connected $s$-decomposable graphs that are not block-decomposable have non-unique decomposition.

\section{Decomposition Rules and Blocks}
 For convenience, we denote an edge that connects nodes $x,y$ by $\overline{xy}$ if the orientation of this edge is unknown or irrelevant, $\overrightarrow{xy}$ if the edge is directed from $x$ to $y$, and $\overleftarrow{xy}$ otherwise.

\begin{dfn}\label{gluingrules}
We recall that a diagram (or graph) is \textit{block-decomposable} (or \textit{decomposable}) if it is obtained by gluing elementary blocks of Table \ref{Decomposition} by the following \textit{gluing rules}:
\begin{enumerate}
\item Two white nodes of two different blocks can be identified. As a result, the graph becomes a union of two parts; the common node is colored black. A white node can neither be identified to itself nor with another node of the same block.
\item A black node can not be identified with any other node.
\item If two white nodes $x$, $y$ of one block (endpoints of edge $\overleftarrow{xy}$) are identified with two white nodes $p$, $q$ of another block (endpoints of edge $\overleftarrow{pq}$), $x$ with $p$, $y$ with $q$ correspondingly, then a multi-edge of weight 2 is formed, and nodes $x=p$, $y=q$ are black.
\item If two white nodes $x$, $y$ of one block (endpoints of edge $\overleftarrow{xy}$) are identified with two white nodes $p$, $q$ of another block (endpoints of edge $\overleftarrow{pq}$), $x$ with $q$, $y$ with $p$ correspondingly, then both edges are removed after gluing, and nodes $x=q$, $y=p$ are black.
\end{enumerate}
\begin{table}[btch]
\centering
\begin{tabular}{ccc}
&&\\
&{\large \textbf{Elementary Blocks}}&{\large \textbf{Triangulation}}\\
\hline
{\large \textbf{Spike:}}&\input{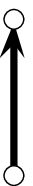}&\input{spikes}\\
{\large \textbf{Triangle:}}&\input{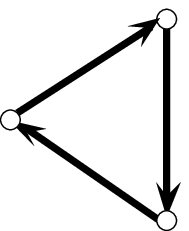}&\input{triangles}\\
{\large \textbf{Infork:}}&\input{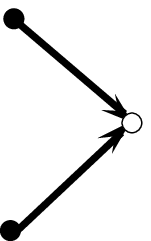}&\input{inforks}\\
{\large \textbf{Outfork:}}&\input{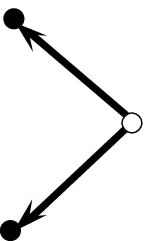}&\input{outforks}\\
{\large \textbf{Diamond:}}&\input{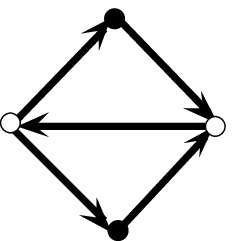}&\input{diamonds}\\
{\large \textbf{Square:}}&\input{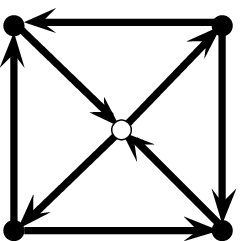}&\input{squares}\\
\hline
\end{tabular}
\caption{Elementary Blocks}\label{Decomposition}
\end{table}
\end{dfn}

\begin{dfn}\label{sdecom}
If a graph $G$ can be obtained by gluing both elementary blocks and new blocks in Table \ref{BOU} by the gluing rules in Definition \ref{gluingrules} and the following new rules, we say the graph is \textit{$s$-decomposable}:
\begin{enumerate}
\item If the graph has multiple edges containing $n$ parallel edges, replace the multiple edge by an edge of weight $2n$. For example, if we glue two parallel spikes of the same direction, we get an edge of weight 4 (see Figure \ref{ME4}).
    \begin{figure}[hctb]
    \centering
    \begin{minipage}[c]{0.5\linewidth}
    \centering
    \includegraphics[width=0.6\linewidth]{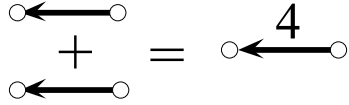}
    \caption{Edge of Weight 4}\label{ME4}
    \end{minipage}
    \end{figure}
\item All single edges have weight 1.
\end{enumerate}
\begin{table}
\begin{center}
\caption{Blocks of Unfolding}\label{BOU}
\begin{tabular}{cccc}
&&&\\
&{\large \textbf{New Blocks}}&{\large \textbf{Unfolding}}&{\large \textbf{Triangulation}}\\
\hline
{\large \textbf{Ia:}}&\input{lws}&\input{uflws}&\input{sflws}\\
{\large \textbf{Ib:}}&\input{rws}&\input{ufrws}&\input{sfrws}\\
{\large \textbf{II:}}&\input{tbt}&\input{uftbt}&\input{sftbt}\\
{\large \textbf{IIIa:}}&\input{lwd}&\input{uflwd}&\input{sflwd}\\
{\large \textbf{IIIb:}}&\input{rwd}&\input{ufrwd}&\input{sfrwd}\\
{\large \textbf{IV:}}&\input{twt}&\input{uftwt}&\input{sftwt}\\
{\large \textbf{V:}}&\input{sq}&\input{ufsq}&\input{sfsq}\\
\hline
\end{tabular}
\end{center}
\end{table}
\end{dfn}

\begin{figure}[hctb]
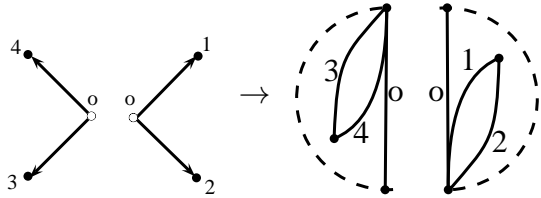

\centering
\input{glue1.tex}\raisebox{3em}{$\rightarrow$}
\input{glue1s.tex}
\caption{Triangulation Gluing}\label{triangglue}
\end{figure}

Gluing two blocks corresponding to gluing two pieces of triangulations of surfaces: gluing two white nodes means gluing the corresponding sides of the triangulations, (see Figure. \ref{triangglue}).

\begin{figure}[hctb]
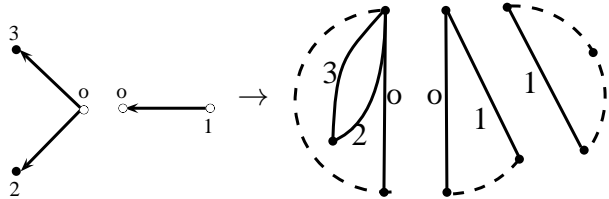

\centering
\input{glue2.tex}\raisebox{3em}{$\rightarrow$}
\input{glue2s.tex}
\caption{Boundary Gluing}\label{boundglue}
\end{figure}

If a decomposable graph has a white node, we will glue a particular piece surface to that node in the corresponding triangulation to form the boundary, see Figure. \ref{boundglue}

It is shown in \cite{FST} that there is a one-to-one correspondence between a decomposition of a graph and an ideal triangulation
of a bordered surfaces with marked points. We show in next section that most graphs with non-unique decomposition correspond to unique bordered surfaces.

\section{Results}
All graphs with non-unique decompositions (s-decompositions) are given in Figure. 78 in \cite{WG}  and Figure. 4 in \cite{WG2}. We list all their block decomposition (s-decomposition) and corresponding ideal (tagged) triangulation of surfaces.

\begin{thm}
If $G$ is a decomposable or $s$-decomposable graph, $G$ is associated to a unique bordered surface unless $G$ is graph 5.
\end{thm}

\begin{table}
\centering
\begin{tabular}{|c|c|c|}
\hline
Graph 1& \multicolumn{2}{|c|}{\includegraphics{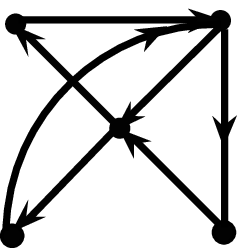}}\\ \hline
Decomposition & \input{1_1.tex} & \input{1_2.tex}\\ \hline
Surfaces &  \input{s1_1.tex} & \input{s1_2.tex}\\
\hline
\end{tabular}
\end{table}

\begin{table}
\centering
\begin{tabular}{|c|c|c|}
\hline
Graph 2& \multicolumn{2}{|c|}{\includegraphics{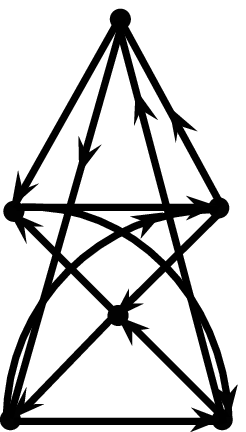}}\\ \hline
Decomposition & \input{2_1.tex} & \input{2_2.tex}\\ \hline
Surfaces &  \input{s2_1.tex} & \input{s2_2.tex}\\
\hline
\end{tabular}
\end{table}

\begin{table}
\centering
\begin{tabular}{|c|c|c|}
\hline
Graph 3& \multicolumn{2}{|c|}{\includegraphics{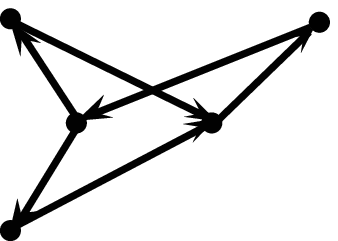}}\\ \hline
Decomposition & \input{3_1.tex} & \input{3_2.tex}\\ \hline
Surfaces &  \input{s3_1.tex} & \input{s3_1.tex} \\
\hline
\end{tabular}
\end{table}
\newpage
\begin{landscape}
\begin{table}
\centering
\begin{tabular}{|c|c|c|c|}
\hline
Graph 4& \multicolumn{3}{|c|}{\input{4.tex}}\\ \hline
Decomposition & \input{4_1.tex} & \input{4_2.tex} & \input{4_3.tex}\\ \hline
Surfaces &  \input{s4_1.tex} &\input{s4_2.tex} & \input{s4_3.tex}\\
\hline
\end{tabular}
\end{table}
\end{landscape}
\begin{table}
\centering
\begin{tabular}{|c|c|c|}
\hline
Graph 5& \multicolumn{2}{|c|}{\includegraphics[width=0.15\linewidth]{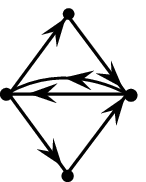}}\\ \hline
Decomposition & \input{5_1.tex} & \input{5_2.tex}\\ \hline
Surfaces &  \input{s5_1.tex} & \input{s5_2.tex} \\
\hline
\end{tabular}
\end{table}
\begin{landscape}
\begin{table}
\centering
\begin{tabular}{|c|c|c|c|}
\hline
Graph 6& \multicolumn{3}{|c|}{\input{6.tex}}\\ \hline
Decomposition & \input{6_1.tex} & \input{6_2.tex} & \input{6_3.tex}\\ \hline
Surfaces &  \input{s6_1.tex} & \input{s6_2.tex} & \input{s6_3.tex} \\
\hline
\end{tabular}
\end{table}

\begin{table}
\centering
\begin{tabular}{|c|c|c|c|}
\hline
Graph 7& \multicolumn{3}{|c|}{\includegraphics[width=0.12\linewidth]{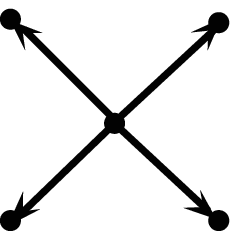}}\\ \hline
Decomposition & \input{7_1.tex} & \input{7_2.tex} &\input{7_3.tex}\\ \hline
Surfaces &  \input{s7_1.tex} &\input{s7_2.tex} &\input{s7_3.tex}\\
\hline
\end{tabular}
\end{table}

\begin{table}
\centering
\begin{tabular}{|c|c|c|c|}
\hline
Graph 7'& \multicolumn{3}{|c|}{\input{7_.tex}}\\ \hline
Decomposition & \input{7_1_.tex} & \input{7_2_.tex} & \input{7_3_.tex} \\ \hline
Surfaces &  \input{s7_1_.tex} & \input{s7_2_.tex} & \input{s7_3_.tex}\\
\hline
\end{tabular}
\end{table}

\begin{table}
\centering
\begin{tabular}{|c|c|c|c|}
\hline
Graph 8& \multicolumn{3}{|c|}{\includegraphics[width=0.12\linewidth]{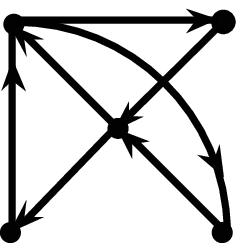}}\\ \hline
Decomposition & \input{8_1.tex} & \input{8_2.tex} & \input{8_3.tex}\\ \hline
Surfaces &  \input{s8_1.tex} & \input{s8_2.tex} & \input{s8_3.tex}\\
\hline
\end{tabular}
\end{table}

\begin{table}
\centering
\begin{tabular}{|c|c|c|c|}
\hline
Graph 9& \multicolumn{3}{|c|}{\includegraphics{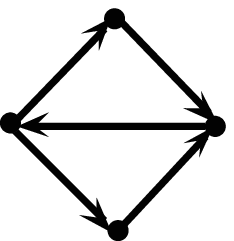}}\\ \hline
Decomposition & \input{9_1.tex} & \input{9_2.tex} & \input{9_3.tex}\\ \hline
Surfaces &  \input{s9_1.tex} & \input{s9_2.tex} &\input{s9_3.tex}\\
\hline
\end{tabular}
\end{table}

\begin{table}
\centering
\begin{tabular}{|c|c|c|c|}
\hline
Graph 10& \multicolumn{3}{|c|}{\includegraphics{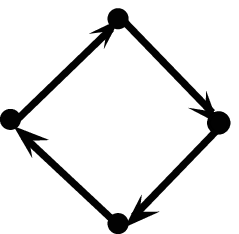}}\\ \hline
Decomposition & \input{10_1.tex} & \input{10_2.tex} & \input{10_3.tex}\\ \hline
Surfaces &  \input{s10_1.tex} & \input{s10_2.tex} &\input{s10_3.tex}\\
\hline
\end{tabular}
\end{table}

\end{landscape}
\begin{table}
\centering
\begin{tabular}{|c|c|c|}
\hline
Graph 11& \multicolumn{2}{|c|}{\includegraphics{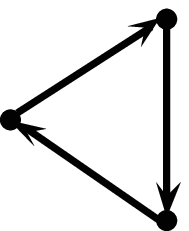}}\\ \hline
Decomposition & \input{11_1.tex} & \input{11_2.tex} \\ \hline
Surfaces &  \input{s11_1.tex} & \input{s11_2.tex} \\
\hline
\end{tabular}
\end{table}

\begin{table}
\centering
\begin{tabular}{|c|c|c|}
\hline
Graph 12& \multicolumn{2}{|c|}{\includegraphics{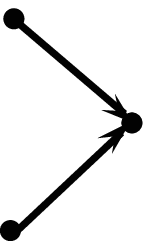}}\\ \hline
Decomposition & \input{12_1.tex} & \input{12_2.tex} \\ \hline
Surfaces &  \input{s12_1.tex} & \input{s12_2.tex} \\
\hline
\end{tabular}
\end{table}

\begin{landscape}
\begin{table}
\centering
\begin{tabular}{|c|c|c|c|}
\hline
Graph 13& \multicolumn{3}{|c|}{\includegraphics{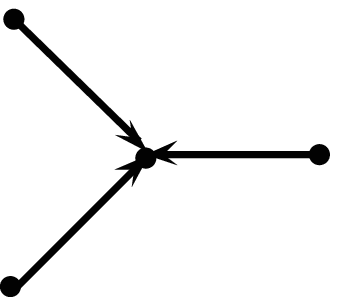}}\\ \hline
Decomposition & \input{13_1.tex} & \input{13_2.tex} & \input{13_3.tex}\\ \hline
Surfaces &  \input{s13_1.tex} & \input{s13_2.tex} &\input{s13_3.tex}\\
\hline
\end{tabular}
\end{table}

\begin{table}
\centering
\begin{tabular}{|c|c|c|c|}
\hline
Graph 14& \multicolumn{3}{|c|}{\includegraphics{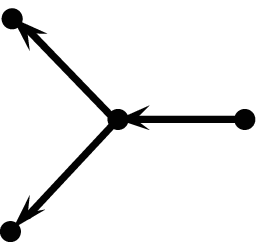}}\\ \hline
Decomposition & \input{14_1.tex} & \input{14_2.tex} & \input{14_3.tex}\\ \hline
Surfaces &  \input{s14_1.tex} & \input{s14_2.tex} &\input{s14_3.tex}\\
\hline
\end{tabular}
\end{table}

\begin{table}
\centering
\begin{tabular}{|c|c|c|c|}
\hline
Graph 15& \multicolumn{3}{|c|}{\input{square.tex}}\\ \hline
Decomposition & \input{square.tex} & \input{15_2.tex} & \input{15_3.tex}\\ \hline
Surfaces &  \input{squares.tex} & \input{s15_2.tex} &\input{s15_3.tex}\\
\hline
\end{tabular}
\end{table}
\end{landscape}
\begin{table}
\centering
\begin{tabular}{|c|c|c|}
\hline
Graph 16& \multicolumn{2}{|c|}{\input{n1.tex}}\\ \hline
Decomposition & \input{n1_1.tex} & \input{n1_2.tex} \\ \hline
Surfaces &  \input{sn1_1.tex} & \input{sn1_2.tex} \\
\hline
\end{tabular}
\end{table}
\begin{table}
\centering
\begin{tabular}{|c|c|c|}
\hline
Graph 17& \multicolumn{2}{|c|}{\input{n2.tex}}\\ \hline
Decomposition & \input{n2_1.tex} & \input{n2_2.tex} \\ \hline
Surfaces &  \input{sn2_1.tex} & \input{sn2_2.tex} \\
\hline
\end{tabular}
\end{table}

\clearpage
\bibliography{mybibs}{}
\bibliographystyle{unsrt}
\end{document}